
\magnification 1200
\input amstex
\documentstyle{amsppt}

\font\bigtenrm=cmr10 scaled\magstep2

\def\non{{{\noindent}}}

\def\ung{\frak{sl}_2}
\def\ch{{\roman{ch}}}
\def\bd{{\bold d}}
\def\Calp{{{{\Cal{P}}}}}
\def\calp{\Calp}
\def\ot{{{\otimes}}}

\def\call{{\Cal{L}}}
\def\superseteq{\supseteq}
\def\non{{{\noindent}}}
\def\l{\lambda}
\def\m{\mu}
\def\n{\nu}

\def\y{{Y(\ung)}}
\def\uqg{\y}
\def\hom{{\roman{Hom}_{\y}}}
\NoBlackBoxes
\nologo
\document
{\bigtenrm{Yangians: their Representations and Characters}}
\vskip24pt\centerline{Vyjayanthi Chari and Andrew Pressley}
\vskip36pt\centerline{\bf Introduction}
\vskip12pt\noindent The Yangian $Y(\frak g)$ associated to a finite-dimensional
complex simple Lie algebra $\frak g$ is a Hopf algebra deformation of the
universal enveloping algebra of the Lie algebra $\frak g[u]$ of polynomial maps
$\Bbb C\to\frak g$ (with Lie bracket defined pointwise). It is important to
understand the finite-dimensional representations of $Y(\frak g)$. For one
thing, they are closely related to rational solutions of the so-called quantum
Yang--Baxter equation (see [7] and [5], Chapter 12). For another, they have
recently been shown to act as \lq quantum symmetry groups\rq\ of certain
integrable systems (see [1]). In [8], Drinfel'd gave a classification of the
finite-dimensional irreducible representations of $Y(\frak g)$ analogous to the
 classification of the finite-dimensional irreducible representations of $\frak
g$ in terms of highest weights, but this left open the problem of giving an \lq
explicit\rq\ realization of these representations. This was!
 solved in [3] when $\frak g=\ung$

. We showed there that every finite-dimensional irreducible representation of
$Y(\ung)$ is a tensor product of representations which are irreducible under
$\ung$, and described the action of $\y$ on these latter representations
completely (for any $\frak g$, $Y(\frak g)$ contains $U(\frak g)$ as a (Hopf)
subalgebra, so any representation of $Y(\frak g)$ can be regarded as a
representation of $\frak g$). Although it is known that, when
${\roman{rank}}(\frak g)>1$, it is not true that every finite-dimensional
irreducible representation of $Y(\frak g)$ is a tensor product of
representations which are irreducible under $\frak g$ (cf. [6]), one can hope
to identify a larger, but still manageable, class of representations of
$Y(\frak g)$ of which every finite-dimensional irreducible representation of
$Y(\ung)$ is a tensor product. Unfortunately, the proof of the tensor product
theorem given in [3] depends on the computation of a certain \lq quantum
Vandermonde determinant\rq, and !
this does not seem to generalise e

asily to the higher rank situation. In the first part of this paper (Section
3), we give a proof of the tensor product theorem in the $\ung$ case using
techniques which, as far as possible, admit generalisation to the rank $>1$
case; in particular, we make no use of the quantum Vandermonde determinant.

A second approach to understanding the finite-dimensional irreducible
representations of $Y(\ung)$ would be to obtain an analogue of the Weyl
character formula for them, i.e. a formula for the character of such a
representation which depends only on its highest weight (the appropriate notion
of character was introduced, and its basic properties obtained, in [9]). In the
second part of this paper (Section 4), we give several such formulas in the
$\ung$ case. Their proof depends of the tensor product theorem described above,
and thus does not generalise to the higher rank case. Nevertheless, one can
hope that the form taken by the formulas in the $\ung$ case might suggest
generalisations for arbitrary $\frak g$, which could then be proved by other
methods.

\vskip36pt\centerline{\bf 1. Yangians}
\vskip12pt\noindent We take the usual basis $\{H,X^+,X^-\}$ of the Lie algebra
$\ung$ (over $\Bbb C$), so that
$$[H,X^\pm]=\pm 2X^\pm,\ \ \ [X^+,X^-]=H.$$
Let $(\ ,\ )$ be the invariant symmetric bilinear form on $\ung$ such that
$$(H,H)=2, \ \ \ (X^+,X^-)=1,$$
and denote by $\Omega$ the Casimir element
$$\Omega=\sum_\l I_\l\ot I_\l,$$
where $\{I_\l\}$ is any orthonormal basis of $\ung$. We also denote by $\Omega$
the element
$$\Omega=\sum_\l I_\l^2$$
in the universal enveloping algebra $U(\ung)$.

\proclaim{Definition 1.1} The Yangian $\y$ is the algebra over $\Bbb C$
generated by elements
$x$, $J(x)$, for $x\in\ung$, with the following
defining relations:
$$\align [x,y]\ \ &\text{(in $\y$)}\ =\ [x,y]\ \ \text{(in $\frak
g$)}\,,\tag1\\ J(ax+by)&=aJ(x)+bJ(y)\,,\tag2\\
[x,J(y)]&=J([x,y])\,,\tag3\\
[J(x),J([y,z])]+&[J(z),J([x,y])]+[J(y),J([z,x])]=\\
&\qquad\sum_{\l,\m,\n}([x,I_\l]\,,\,
[[y,I_\m],[z,I_\n]])
\{I_\l,I_\m,I_\n\}\,,\tag4\\
[[J(x),J(y)]\,,\,&[z,J(w)]]+[[J(z),J(w)]\,,\,[x,J(y)]]=\\
& \sum_{\l,\m,\n}([x,I_\l]\,,\,[[y,I_\m],
[[z,w],I_\n]])\{I_\l,I_\m,J(I_\n)\}\,,\tag5
\endalign$$ for all $x$, $y$, $z\in\ung$, $a$,
$b\in\Bbb C$. Here, for any elements $z_1$, $z_2$, $z_3\in\y$, we set
$$\{z_1,z_2,z_3\}=\frac1{24}\sum_\pi
z_{\pi(1)}z_{\pi(2)}z_{\pi(3)},$$ the sum being over all
permutations $\pi$ of $\{1,2,3\}$.

The Yangian $\y$ has a Hopf algebra structure with counit $\epsilon$,
comultiplication $\Delta$ and antipode $S$ given by
$$\align
\Delta(x)&=x\otimes 1+1\otimes x,\tag6\\
\Delta(J(x))&=J(x)\otimes 1+1\otimes J(x)+\frac12[x\otimes
1\,,\,\Omega],\tag7\\
S(x)&=-x,\ \ S(J(x))=-J(x)+ x,\tag8\\
\epsilon(x)&=\epsilon(J(x))=0.\tag9\endalign$$
\endproclaim

We shall also need the following presentation of $\y$, given in [8]:
\proclaim{Theorem 1.2} The
Yangian $\y$ is isomorphic to the
associative algebra with generators $X_{k}^{{}\pm{}}$,
$H_{k}$, $k\in\Bbb N$, and the following
defining relations:
$$\align
[H_k,H_l]=0,\ \ [H_0,X_k^\pm]&=\pm 2X_k^\pm,\ \ [X_k^+,X_l^-]=H_{k+l},\tag10\\
[H_{k+1},X_l^\pm]-[H_k,X_{l+1}^\pm]&=\pm(H_kX_l^\pm+X_l^\pm H_k),\tag11\\
[X_{k+1}^\pm,X_l^\pm]-[X_k^\pm,X_{l+1}^\pm]&=\pm(X_k^\pm X_l^\pm+X_l^\pm
X_k^\pm).\tag12\endalign$$
The isomorphism $f$ between the two realizations of $\y$ is given by
$$\align
f(H)&=H_0,\ \  f(X^\pm)=X_0^\pm,\\
f(J(H))&=H_1+\frac12(X_0^+X_0^-+X_0^-X_0^+-H_0^2),\\
f(J(X^\pm))&=X_1^\pm-\frac14(X_0^\pm H_0+H_0X_0^\pm).\qed\endalign$$
\endproclaim

The presentation 1.1 of $\y$ shows that there is a canonical map $\ung\to\y$
(it is known that this map is injective). Thus, any $\y$-module may be regarded
as an $\ung$-module.

We shall make use of the following automorphism of $\y$.
\proclaim{Proposition 1.3} There is a one-parameter group $\{\tau_a\}_{a\in\Bbb
C}$ of Hopf algebra automorphisms of $\y$ given in terms of the presentation
1.1 by
$$\tau_a(x)=x,\ \ \tau_a(J(x))=J(x)+ax,$$
for $x\in\ung$, and in terms of the presentation 1.2 by
$$\tau_a(H_{k})=\sum_{r=0}^k\left({k\atop r}\right)a^{k-r}H_{r},\ \ \
\tau_a(X_{k}^\pm)=\sum_{r=0}^k\left({k\atop r}\right)a^{k-r}X_{r}^\pm.\ \qed$$
\endproclaim

This is easily checked using 1.1 and 1.2.

We shall also need the following weak version of the Poincar\'e--Birkhoff--Witt
theorem for $\y$.
\proclaim{Proposition 1.4} Let $Y^+$, $Y^-$ and $Y^0$ be the subalgebras of
$\y$ generated by the $X_{k}^+$, the $X_{k}^-$ and the $H_{k}$, respectively
($k\in \Bbb N$). Then,
$$\y=Y^-.Y^0.Y^+.\ \ \ \qed$$
\endproclaim

The proof is straightforward.

\vskip36pt\centerline{\bf 2. Representations}
\vskip 12pt\non
If $W$ is an $\ung$-module, a non-zero vector $w\in W$ is said to be of weight
$r$ if $H.w=rw$, and is said to be an $\ung$-highest weight vector if, in
addition, $X^+.w=0$; if $W=U(\ung).w$, then $W$ is called a highest weight
$\ung$-module with highest weight $r$. Lowest weight vectors and $\ung$-modules
are defined similarly. For any $r\in\Bbb N$, denote by $W_r$ the unique
irreducible highest weight $\ung$-module with highest weight $r$.

\vskip6pt
Suppose now that $V$ is a $\y$-module. A non-zero vector $v\in V$ is called a
$\y$-highest weight vector if $v$ is an eigenvector of $H_k$, say
$$H_k.v=d_kv,\ \ \ \ (d_k\in\Bbb C)$$
and is annihilated by $X_k^+$, for all $k\in\Bbb N$. Then, $V$ is called a
highest weight $\y$-module if $V=\y.v$ for some $\y$-highest weight vector
$v\in V$, and the $\Bbb N$-tuple of scalars $\bold{d}=\{d_{k}\}_{k\in\Bbb N}$
is called its highest weight. It is not difficult to show that, for every
${\bold d}\in\Bbb C^{\Bbb N}$, there is an irreducible $\y$-module $V({\bold
d})$, unique up to isomorphism, such that $V({\bold d})$ has highest weight
${\bold d}$. Lowest weight vectors and modules for $\y$ are defined similarly.

The following theorem of Drinfel'd [8] classifies the finite-dimensional
irreducible $\y$-modules. Let $\Calp$ be the set of monic polynomials in $\Bbb
C[u]$, where $u$ is an indeterminate.
\proclaim{Theorem 2.1} (i) Every finite-dimensional irreducible $\y$-module is
both highest weight and lowest weight.

(ii) If ${\bold d}=\{d_k\}_{k\in\Bbb N}\in\Bbb C^{\Bbb N}$, the $\y$-module
$V({\bold d})$ is finite-dimensional if and only if there exists $P\in\Calp$
such that
$$\frac{P(u+1)}{P(u)}=1+\sum_{k=0}^{\infty}d_{k}u^{-k-1},\tag13$$
in the sense that the right-hand side is the Laurent expansion of the left-hand
side about $u=\infty$. $\qed$
\endproclaim

If $V$ is a finite-dimensional irreducible $\uqg$-module, we call the
associated polynomial $P$ the Drinfel'd polynomial of $V$.

More generally, if $V$ is any finite-dimensional $\uqg$-module and $v\in V$ is
a $\y$-highest weight vector, with
$$H_{k}.v=d_{k} v$$
for some $d_{k}\in\Bbb C$ and all $k\in\Bbb N$, it follows from 2.1 that there
exists $P\in\Calp$ such that
$$\frac{P(u+1)}{P(u)} =1+\sum_{k=0}^{\infty}d_{k}u^{-k-1}.$$

\proclaim{Proposition 2.2} Let $V_1$, $V_2$ be finite-dimensional
$\uqg$-modules, and let $v_1\in V_1$, $v_2\in V_2$ be $\y$-highest weight
vectors with associated polynomials $P_1$ and $P_2$. Then, $v_1\ot v_2$ is a
$\y$-highest weight vector in $V_1\ot V_2$ with associated polynomial $P_1P_2$.
\ \ \qed\endproclaim

This is Proposition 4.6 in [3].

Given a finite-dimensional $\y$-module $V$, we can define the following
associated $\y$-modules:

(i) $V(a)$: this is obtained pulling back $V$ through $\tau_a$;

(ii) the left dual $^tV$ and right dual $V^t$: these are given by the following
actions of $\y$ on the vector space dual of $V$:
$$\align (y.f)(v)&=f(S(y).v),\ \ \ \ \ y\in\uqg ,f\in {}^tV, v\in V,\\
(y.f)(v)&=f(S^{-1}(y)).v,\ \ \ \ \ y\in\uqg ,f\in V^t, v\in V.\endalign$$
Clearly, if $V$ is irreducible, so are all the representations defined above.

\proclaim{Proposition 2.3} Let $U$, $V$ and $W$ be finite-dimensional
$\y$-modules, and let $a\in\Bbb C$. Then,
\vskip6pt
\non(i) ${\roman{Hom}}_{\y}(U,V\ot W)\cong {\roman{Hom}}_{\y}(^t V\ot U, W)$;

\non(ii)  ${\roman{Hom}}_{\y}(U,W\ot V)\cong {\roman{Hom}}_{\y}(U\ot V^t, W)$.
\qed\endproclaim

The proof is straightforward, using the fact that $S$ is a coalgebra
anti-automorphism of $\y$.

The following result describes the Drinfel'd polynomials of the modules defined
above. See [4] for the proof.

\proclaim{Proposition 2.4}  Let $V$ be a finite-dimensional irreducible
$\uqg$-module with Drinfel'd polynomial $P$, and let $a\in\Bbb C$.
Then:
\vskip6pt\non
(i) The Drinfel'd polynomial of $V(a)$ is
$P(u-a)$.

\non(ii) The Drinfel'd polynomials of $^tV$ and $V^t$ are
$P(u+1)$ and $P(u-1)$, respectively. \qed
\endproclaim

\proclaim{Corollary 2.5} ${}^t V\cong V(-1)$ and $V^t\cong
V(1)$.\qed\endproclaim

We conclude this section with the following result.
\proclaim{Proposition 2.6} Let $V$ be a finite-dimensional $\uqg$-module. Then,
$V$ is irreducible if and only if $V$ and $^tV$ (resp. $V$ and $V^t$) are both
highest weight $\uqg$-modules.
\endproclaim
\demo{Proof} The \lq only if\rq\  part follows from 2.1 (i). For the converse,
suppose that $V$ and $^tV$ are highest weight (the other case is identical).
Let $v$ be a $\y$-highest weight vector in $V$ of weight $n$ (say) for $\ung$.
Let $0\ne W$ be an irreducible $\y$-submodule of $V$, and let $m$ (say) be the
highest weight of $W$ as an $\ung$-module; thus, $m\le n$. Then, $^tW$ is a
quotient of $^tV$, and these $\ung$-modules have maximal weights $m$ and $n$,
respectively. Since $^tV$ is a $\y$-highest weight module, its highest weight
vector must map to a non-zero element of $^tW$. Hence, $n\le m$. Thus, $m=n$
and $W=V$.\qed\enddemo

\vfill\eject\centerline{\bf 3. Classification}
\vskip12pt\noindent For any $a\in\Bbb C$, $r\in\Bbb N$, the {\it string}
$S_r(a)$ is the set of complex numbers
$$S_r(a)=\{a, a+1,a+2,\ldots,a+r-1\}.$$
(Note that this is different from the notation used in [3].) We say that
$S_r(a)$ begins at $a$ and ends at $a+r-1$; its cardinality $|S_r(a)|=r$ is
often called the length of $S_r(a)$.

Two strings $S$ and $T$ are said to be in special position if $S\cup T$ is a
string which is strictly longer than both $S$ and $T$; otherwise, $S$ and $T$
are in general position. We leave it to the reader to prove

\proclaim{Lemma 3.1} Let $a$, $b\in\Bbb C$, $p$, $q\in\Bbb N$. Then, $S_p(a)$
and $S_q(b)$ are in general position if and only if
\vskip6pt\non(i) $S_p(a)\subseteq S_q(b)$ or $S_q(b)\subseteq S_p(a)$, or

\non(ii) $S_p(a)\cap S_q(b)=S_p(a+1)\cap S_q(b)=S_p(a-1)\cap S_q(b)=\emptyset$.
\qed\endproclaim

The following result contains the combinatorial properties of strings that we
shall need. Again, we leave the proof to the reader.

\proclaim{Proposition 3.2} Let $S_1,S_2,\ldots,S_p$ ($p\ge 1$) be a collection
of strings, every pair of which are in general position.
\vskip6pt\non(i) Let $a\in\Bbb C$ be such that, for all $n\ge 0$,
$a+n\notin\bigcup_{q=1}^p S_q$. Then,
\vskip6pt(a) if $a-1\notin \bigcup_{q=1}^p S_q$, the string $\{a\}$ is in
general position with respect to every $S_q$;

(b) if $a-1\in S_r$, and if $|S_r|$ is maximal with this property, then
$S_r\cup \{a\}$ is a string in general position with respect to every $S_q$.
\vskip6pt\non(ii) Let $a\in S_{t}$, and assume that $|S_t|$ is minimal with
this property. Then, $S_t\backslash\{a\}$ is a string in general position with
respect to every $S_q$. \qed\endproclaim

A {\it multiset} is a map ${\Cal S}\to \Bbb N$, where ${\Cal S}$ is a finite
set of complex numbers. Define the union, intersection and cardinality of
multisets in the obvious way.

\proclaim{Corollary 3.3} Any multiset $S$ can be written uniquely as a union of
strings, any two of which are in general position. \endproclaim
\demo{Proof} By induction on the cardinality of $S$, using 3.2
(ii).\qed\enddemo

We call the decomposition of $S$ into strings given by 3.3 its canonical
decomposition.

\vskip12pt
The following $\y$-modules are the \lq building blocks\rq\ out of which an
arbitrary finite-dimensional $\y$-module will be constructed by taking tensor
products.

\proclaim{Definition 3.4} If $a\in\Bbb C$, $r\in\Bbb N$, then $W_r(a)$ is the
finite-dimensional irreducible $\y$-module whose Drinfel'd polynomial is
$$P_{r,a}(u)=\prod_{b\in S_r(a)}(u-b).$$
\endproclaim

The explicit action of the generators $H_k$, $X_k^\pm$ ($k\in\Bbb N$) on a
suitable basis of $W_r(a)$ is given in Proposition 2.6 and Corollary 2.7 in
[3]. We recall the result:

\proclaim{Proposition 3.5} For any $r\ge 1$, $a\in\Bbb C$, $W_r(a)$ has a basis
$\{w_0,w_1,\ldots,w_r\}$ on which the action of $\y$ is given by
$$\align
X_k^+.w_s&=(s+a)^k(s+1)w_{s+1},\ \ X_k^-.w_s=(s+a-1)^k(r-s+1)w_{s-1},\\
H_k.w_s&=\left((s+a-1)^ks(r-s+1)-(s+a)^k(s+1)(r-s)\right)w_s\endalign$$
(it is understood that $w_{-1}=w_{r+1}=0$). In particular,  $W_r(a)$ is
isomorphic to $W_r$ as an $\ung$-module. \qed\endproclaim

Using these formulas, and the formula for the comultiplication of $\y$ given in
1.1, it is straightforward to prove

\proclaim{Proposition 3.6} Let $a\in\Bbb C$ and let $v_0=v^+\ot v^- -v^-\ot
v^+$, where $v^+$ (resp. $v^-$) is an $\ung$-highest (resp. lowest) weight
vector in $W_1$.
\vskip6pt\non(i) $W_1(a+1)\ot W_1(a)=\y.(v^+\ot v^-)$, and we have a short
exact sequence of $\y$-modules
$$0\to\y.v_0\to W_1(a+1)\ot W_1(a)\to W_2(a)\to 0,$$
where $\y.v_0$ is the one-dimensional trivial module.

\non(ii) $W_1(a)\ot W_1(a+1)=\y.v_0$ and we have a short exact sequence of
$\y$-modules
$$0\to \y.(v^+\ot v^-)\to W_1(a)\ot W_1(a+1)\to\Bbb C\to 0,$$
where $\y.(v^+\ot v^-)\cong W_2(a)$.

\non(iii) ${}^tW_1(a)\cong W_1(a-1)$, $W_1(a)^t\cong W_1(a+1)$.\qed\endproclaim

We use these computations to prove

\proclaim{Proposition 3.7} Let $a_1,a_2,\ldots,a_r\in\Bbb C$, $r\ge 1$. Then:
\vskip6pt\non(i) if $a_j-a_i\ne 1$ when $i<j$,
$$W_1(a_1)\ot W_1(a_2)\ot\cdots\ot W_1(a_r)$$
is a highest weight $\y$-module;

\non(ii) if $a_i-a_j\ne 1$ when $i<j$,
$$W_1(a_1)\ot W_1(a_2)\ot\cdots\ot W_1(a_r)$$
does not contain a $\y$-highest weight vector of weight $<r$ for $\ung$.
\endproclaim
\demo{Proof} (i) By induction on $r$. If $r=1$, there is nothing to prove, and
the $r=2$ case is contained in 3.6. Assume now that $r>2$ and that the result
is known for $r-1$. If $W_1(a_1)\ot W_1(a_2)\ot\cdots\ot W_1(a_r)$ is not a
highest weight $\y$-module, it has an irreducible quotient $V(P)$, say, so that
$$\hom(W_1(a_1)\ot W_1(a_2)\ot\cdots\ot W_1(a_r),V(P))\ne 0.$$
By 2.3 and 3.6,
$$\hom(W_1(a_2)\ot\cdots\ot W_1(a_r), W_1(a_1+1)\ot V(P))\ne 0.$$
Let $F$ be a non-zero element of this space of homomorphisms. By the induction
hypothesis,
$$W_1(a_2)\ot\cdots\ot W_1(a_r)=\y.(v^+)^{\ot r-1},$$
where $v^+$ denotes an $\ung$-highest weight vector in $W_1$, so $F((v^+)^{\ot
r-1})$ must be a non-zero multiple of $v^+\ot v_P$, where $v_P$ is a
$\y$-highest weight vector in $V(P)$. By 2.2, $a_1+1=a_i$ for some $i\ge 2$,
contradicting our assumption.
\vskip6pt(ii) Let $0\ne v\in W_1(a_1)\ot W_1(a_2)\ot\cdots\ot W_1(a_r)$ be a
$\y$-highest weight vector of weight $s$ for $\ung$, and let $V=\y.v$. Then,
$$\hom(V,W_1(a_1)\ot W_1(a_2)\ot\cdots\ot W_1(a_r))\ne 0,$$
and so, by 2.3 and 3.6,
$$\hom(W_1(a_r-1)\ot\cdots\ot W_1(a_1-1),{}^t V)\ne 0.$$
By part (i), $W_1(a_r-1)\ot W_1(a_2-1)\ot\cdots\ot W_1(a_1-1)$ is a
$\y$-highest weight module, so ${}^t V$ contains a vector of weight $r$ for
$\ung$. It follows that $r\le s$, and hence that $r=s$. \qed\enddemo

This result has several consequences.

\proclaim{Corollary 3.8} Let $P_1,P_2,\ldots P_r\in\calp$, $r\ge 1$, and assume
that, if $a_i$ is a root of $P_i$ and $a_j$ a root of $P_j$, where $i<j$, then
$a_j-a_i\ne 1$. Then,
$$V(P_1)\ot V(P_2)\ot\cdots\ot V(P_r)$$
is a highest weight $\y$-module.
\endproclaim
\demo{Proof} Let $\{a_{1i},a_{2i},\ldots,a_{d_ii}\}$ be the multiset of roots
of $P_i$, where $d_i={\roman{deg}}(P_i)$ and each root is repeated according to
its multiplicity. Order the roots so that $a_{ti}-a_{si}\ne 1$ if $s<t$.
Then, by 2.2, $V(P_s)$ is a quotient of
$$V_i=W_1(a_{1i})\ot W_1(a_{2i})\ot\cdots\ot W_1(a_{d_ii}),$$
and by 3.7 (i),
$$V_1\ot V_2\ot\cdots\ot V_r$$
is a highest weight $\y$-module. Hence, $V(P_1)\ot\cdots\ot V(P_r)$ is a
highest weight $\y$-module. \qed\enddemo

\proclaim{Corollary 3.9} Every finite-dimensional irreducible $\y$-module is
(isomorphic to) a quotient (resp. a submodule) of a tensor product of modules
of the form $W_1(a)$, for $a\in\Bbb C$. \endproclaim
\demo{Proof} If $P\in\calp$, let $d={\roman{deg}}(P)$ and order the roots
$a_1,a_2,\ldots,a_d$ of $P$ (repeated according to multiplicity) so that
$a_j-a_i\ne 1$ if $i<j$. Then, $V(P)$ is a quotient of
$$W_1(a_1)\ot W_1(a_2)\ot\cdots\ot W_1(a_d).$$

The corresponding statement about submodules follows by taking (left or right)
duals. \qed
\enddemo

To obtain the final consequence of 3.7, we need

\proclaim{Lemma 3.10} Let $V$ and $W$ be finite-dimensional irreducible
$\y$-modules. Then, $V\ot W$ is irreducible if and only if $V\ot W$ and $W\ot
V$ are both highest weight $\y$-modules. In particular, $V\ot W$ is irreducible
if and only if $W\ot V$ is irreducible.\endproclaim
\demo{Proof} Let $v^+$ and $w^+$ be $\y$-highest weight vectors in $V$ and $W$,
respectively, and let their Drinfel'd polynomials be $P$ and $Q$. Assume that
$V\ot W$ is irreducible. The irreducible quotient of the submodule of $W\ot V$
generated by $w^+\ot v^+$ has Drinfel'd polynomial $PQ$, by 2.2, and hence is
isomorphic to $V\ot W$, by 2.1. For dimensional reasons, the subquotient must
therefore be $W\ot V$. Hence, $W\ot V$ is irreducible. That both $V\ot W$ and
$W\ot V$ are highest weight now follows from 2.1.

Conversely, assume that $V\ot W$ and $W\ot V$ are highest weight. If $V\ot W$
is reducible, it contains an irreducible $\y$-submodule $Z$, say, whose maximal
weight as an $\ung$-module is strictly less than that of $V\ot W$. Using 2.3,
we get a non-zero homomorphism
$$W^t\ot V^t\to Z^t.$$
Using 2.5 and twisting by $\tau_{-1}$, we get a non-zero homomorphism
$$W\ot V\to Z.$$
This contradicts the fact that $W\ot V$ is highest weight. \qed\enddemo

The following result is now immediate from 2.3, 3.8 and 3.10.

\proclaim{Corollary 3.11} Let $P_1,\ldots,P_r\in\calp$, $r\ge 1$, and assume
that, if $a_i$ is a root of $P_i$, $a_j$ a root of $P_j$, and $i<j$, then $1\ne
a_i-a_j\ne -1$. Then,
$$V(P_1)\ot\cdots\ot V(P_r)$$
is an irreducible $\y$-module. \qed\endproclaim

We are now in a position to take the crucial step towards the classification of
the finite-dimensional irreducible $\y$-modules.

\proclaim{Proposition 3.12} Let $a_1,\ldots,a_p\in\Bbb C$,
$r_1,\ldots,r_p\in\Bbb N$, $p\ge 1$, and assume that $S_{r_i}(a_i)\subseteq
S_{r_1}(a_1)$ for all $i=1,\ldots,p$. Then,
$$W_{r_1}(a_1)\ot\cdots\ot W_{r_p}(a_p)$$
is an irreducible $\y$-module. \endproclaim
\demo{Proof} By induction on $p$. If $p=1$, there is nothing to prove. Now
assume that $p=2$. To simplify the notation in this case, we consider
$W_k(a)\ot W_l(b)$ instead of $W_{r_1}(a_1)\ot W_{r_2}(a_2)$, and assume that
$S_l(b)\subseteq S_k(a)$ (so $l\le k$).

Suppose first that $l=1$. If $k=1$, then $a=b$ and the result follows from
3.11. If $k=2$, then $b=a$ or $a+1$. Assume that $b=a+1$ (the other case is
similar). By 3.8, $W_1(a+1)\ot W_2(a)$ is a highest weight $\y$-module, since
the roots of the Drinfel'd polynomial of $W_2(a)$ are $a$ and $a+1$. By 3.10,
it suffices to prove that $W_2(a)\ot W_1(a+1)$ is a highest weight $\y$-module.
Assuming otherwise, $W_2(a)\ot W_1(a+1)$ has an irreducible quotient
$\y$-module, which must have highest weight 1 as an $\ung$-module, and hence
must be isomorphic as a $\y$-module to $W_1(c)$, for some $c\in\Bbb C$. By 2.3
and 2.5, this implies the existence of a non-zero homomorphism of $\y$-modules
$$F:W_2(a)\to W_1(c)\ot W_1(a).$$
Since $F$ must, for weight reasons, map the highest weight vector in $W_2(a)$
to a non-zero multiple of the tensor product of the highest weight vectors in
$W_1(c)$ and $W_1(a)$, 2.2 implies that $c=a+1$, contradicting 3.6 (i).

Assume now that $k>2$ and that the result is known for smaller values of $k$
(we are still assuming that $l=1$). We consider three cases:
\vskip6pt\non{\it Case I: $b\ne a$ and $b\ne a+k-1$.} Then,
$$S_{k-1}(a+1)\superseteq S_1(b)\subseteq S_{k-1}(a),$$
so by the induction hypothesis on $k$,
$$W_{k-1}(a)\ot W_1(b)\ \ \text{and}\ \ W_1(b)\ot W_{k-1}(a+1)$$
are both irreducible $\y$-modules. By 3.8 and the assumption $b\ne a+k-1$,
$$W_1(a+k-1)\ot W_{k-1}(a)\ot W_1(b)$$
is a highest weight $\y$-module, and hence so is its quotient $W_k(a)\ot
W_1(b)$ (note that $W_1(a+k-1)\ot W_{k-1}(a)$ is highest weight by 3.8).
Similarly, by considering
$$W_1(b)\ot W_{k-1}(a+1)\ot W_1(a)$$
and using the assumption $b\ne a$, one sees that
$W_1(b)\ot W_k(a)$ is highest weight. Lemma 3.10 completes the proof.
\vskip6pt\non{\it Case II: $b=a+k-1$}. By the induction hypothesis on $k$,
$W_1(b)\ot W_{k-1}(a+1)$ is irreducible. On the other hand, since $k>2$, by 3.7
and 3.10,
$$W_1(b)\ot W_1(a)\cong W_1(a)\ot W_1(b)$$
as $\y$-modules, so
$$W_{k-1}(a+1)\ot W_1(a)\ot W_1(b)\cong
W_{k-1}(a+1)\ot W_1(b)\ot W_1(a)\tag14$$
as $\y$-modules. By 3.8, the right-hand side of (14) is highest weight, hence
so is the left-hand side. Hence, its quotient $W_k(a)\ot W_1(b)$ is highest
weight. That $W_1(b)\ot W_k(a)$ is highest weight is immediate from 3.8, so
3.10 again completes the proof.
\vskip6pt\non{\it Case III: $b=a$}. This is similar to Case II. We omit the
details.
\vskip6pt
We have now proved the result when $l=1$ (and $p=2$). We next assume that $l>1$
and that the result is known for smaller values of $l$. We prove the result for
$l$ by induction on $k$, starting at $k=l$.

If $k=l$, then $a=b$ and the induction hypothesis on $l$ gives that
$$W_1(a+l-1)\ot W_{l-1}(a)\ot W_l(a)$$
is highest weight, and hence so is its quotient $W_l(a)\ot W_l(a)$. By 3.10,
this last module is irreducible.

For the inductive step, we distinguish three cases:
\vskip6pt\non{\it Case I: $b-a\ne 0$ or $k-l$}. Then,
$$S_{k-1}(a+1)\superseteq S_l(b)\subseteq S_{k-1}(a).$$
By the induction hypothesis on $k$,
$$W_l(b)\ot W_{k-1}(a+1)\ \ \text{and}\ \ W_{k-1}(a)\ot W_l(b)$$
are both irreducible. By 3.8,
$$W_l(b)\ot W_{k-1}(a+1)\ot W_1(a)\ \ \text{and}\ \
W_1(a+1)\ot W_{k-1}(a)\ot W_l(b)$$
are both highest weight, and hence so are their quotients
$$W_l(b)\ot W_k(a)\ \ \text{and}\ \ W_k(a)\ot W_l(b).$$
{\it Case II: $b-a=k-l$}. By the induction hypothesis on $l$, $W_{l-1}(b)\ot
W_k(a)$ is irreducible. By 3.8,
$$W_1(b+l-1)\ot W_{l-1}(b)\ot W_k(a)$$
is highest weight, hence so is $W_l(b)\ot W_k(a)$. On the other hand,
$$\aligned
W_1(b+l-1)\ot W_{l-1}(b)\ot W_k(a)&\cong
W_1(b+l-1)\ot W_k(a)\ot W_{l-1}(b)\\
&\cong W_k(a)\ot W_1(b+l-1)\ot W_{l-1}(b)
\endaligned\tag15$$
as $\y$-modules: the first isomorphism uses the fact that $S_{l-1}(b)$ and
$S_k(a)$ are in general position, the second that $\{b+l-1\}$ and $S_k(a)$ are
in general position (and both isomorphisms use the induction hypothesis on
$l$). But we saw above that the first tensor product in (15) is highest weight,
hence so is the last, and hence so is its quotient $W_k(a)\ot W_l(b)$.
\vskip6pt\non{\it Case III: $b=a$}. This is similar to Case II.
\vskip6pt
We have now proved 3.12 in the case $p=2$. Assume next that $p>2$ and that the
result is known for smaller values of $p$. Let $S$ be the union of the strings
$S_{r_i}(a_i)$, $i=1,\ldots,p$, considered as a set with multiplicities. We
prove the result for $p$ by induction on $|S|$. If $|S|=p$, we are considering
a tensor product of the form $W_1(a)^{\ot p}$, which is irreducible by 3.11.

Assume now that $|S|>p$ and that the result is known for smaller values of
$|S|$. Note that $a_1+r_1-1\in S$ but $a_1+n\notin S$ if $n\ge r_1$. Let
$S'=S\backslash\{a_1+r_1-1\}$, and let $i$ be such that $a_1+r_1-1\in
S_{r_i}(a_i)$ and such that $r_i$ is minimal with this property. By 3.2 (ii),
$$S'=\bigcup_{j\ne
i}S_{r_j}(a_j)\cup\left(S_{r_i}(a_i)\backslash\{a_1+r_1-1\}\right)$$
is the canonical decomposition of $S'$. By the induction hypothesis on $|S|$,
$$W_{r_i-1}(a_i)\ot\bigotimes_{j\ne i}W_{r_j}(a_j)\tag16$$
is irreducible. (Note that, by 3.10 and the induction hypothesis on $p$, the
second tensor product in (16) is independent of the order of the factors, up to
isomorphism.) By 3.8,
$$W_1(a_1+r_1-1)\ot W_{r_i-1}(a_i)\ot\bigotimes_{j\ne i}W_{r_j}(a_j)$$
is highest weight, hence so is its quotient
$$W_{r_i}(a_i)\ot\bigotimes_{j\ne i} W_{r_j}(a_j).$$
By the $p=2$ case, this module is unchanged, up to isomorphism, by successively
interchanging adjacent factors in the tensor product, so we deduce that
$$\bigotimes_{j\ne i} W_{r_j}(a_j)\ot W_{r_i}(a_i)$$
is also highest weight. The usual application of 3.10 completes the proof.
\qed\enddemo

We can now prove the main result of this section.

\proclaim{Theorem 3.13} Let $P\in\calp$, and write the multiset $S(P)$ of roots
of $P$ as a union of strings in general position, say
$$S(P)=\bigcup_{i=1}^p S_{r_i}(a_i).$$
Then, as $\y$-modules,
$$V(P)\cong \bigotimes_{i=1}^p W_{r_i}(a_i)$$
(the factors in the tensor product can be taken in any order).\endproclaim
\demo{Proof} By induction on $p$. If $p=1$, there is nothing to prove. Assume
now that $p>1$ and that the result is known for smaller values of $p$.

If all the strings $S_{r_i}(a_i)$, $i=1,\ldots,p$, are contained in a single
string, the result was proved in 3.12. Otherwise, let $r_i$ be maximal among
$r_1,\ldots,r_p$, and let
$$S'=\bigcup_{\{j|S_{r_j}(a_j)\subseteq S_{r_i}(a_i)\}}S_{r_j}(a_j),\ \ \
S''=S\backslash S'$$
(with $S$, $S'$ and $S''$ considered as sets with multiplicities). Then,
$S=S'\cup S''$ and
$$S'\cap S''=(S'+1)\cap S''=(S'-1)\cap S''=\emptyset.$$
Let $P'$, $P''\in\calp$ have multisets of roots $S'$ and $S''$, respectively.
By 3.11, $V(P')\ot V(P'')$ is irreducible. By the induction hypothesis, $V(P')$
and $V(P'')$ are both isomorphic to tensor products as in the statement of the
theorem, so an application of 3.10 to re-order the factors, if necessary,
completes the proof. \qed\enddemo
\vskip9pt\noindent{\it Remark} It is instructive to consider the classical
analogue of 3.13. As we mentioned in the Introduction, the classical analogue
of $\y$ is $U(\ung[u])$.

\proclaim{Theorem 3.14} Every finite-dimensional irreducible $\ung[u]$-module
$V$ is generated by a vector $v$ such that
$$(H\ot u^k).v=d_k v,\ \ \ \ (X^+\ot u^k).v=0$$
for some $d_k\in\Bbb C$ and all $k\in\Bbb N$. Moreover, there exists a monic
polynomial $P\in\Bbb C[u]$ such that
$$\sum_{k=0}^\infty d_ku^{-k-1}=\frac1{P(u)}\frac{dP}{du}. \qed$$
\endproclaim

It is not difficult to prove (cf. [2]) that every finite-dimensional
irreducible $\ung[u]$-module is a tensor product of modules of the form
$W_r(a)$, for some $r\ge 1$, $a\in\Bbb C$, where $W_r(a)$ is obtained by
pulling back $W_r$ via the Lie algebra homomorphism $\ung[u]\to\ung$ given by
setting $u=a$. The polynomial associated to $W_r(a)$ is $(u-a)^r$, and the
polynomial is multiplicative on irreducible tensor products (cf. 2.2).

\vskip36pt\centerline{\bf 4. Characters}
\vskip12pt\noindent The appropriate definition of the character of a
finite-dimensional $\y$-module was given in [9]. Let $\call$ be the subgroup of
the group of units of the ring $\Bbb C[[u^{-1}]]$ of the form
$$f_\bd(u)=1+\sum_{k=0}^\infty d_ku^{-k-1},$$
for some $\bd=\{d_k\}_{\in\Bbb C}\in\Bbb C^{\Bbb N}$,
and let $\Bbb C[\call]$ be its group algebra. If $f\in\call$, we write $e(f)$
for the corresponding basis element of $\Bbb C[\call]$.

If $V$ is a finite-dimensional $\y$-module, and $\bold d=\{d_k\}_{k\in\Bbb N}$,
set
$$
V_{\bold d}=\{v\in V  \mid(H_k-d_k)^N.v=0\ \text{for $N>>0$}\},$$
and if $P\in\calp$ is such that
$$\frac{P(u+1)}{P(u)}=f_\bd(u)$$
(cf. 2.1 (ii)), set $e(P)=e(f_\bd)$ (abusing notation a little).

\proclaim{Definition 4.1} If $V$ is a finite-dimensional $\y$-module, its
character is
$$\ch(V)=\sum_\bd{\roman{dim}}(V_\bd)e(f_\bd).$$\endproclaim

The main result proved in [9] is

\proclaim{Proposition 4.2} (i) If
$$0\to U\to V\to W\to 0$$
is a short exact sequence of $\y$-modules, then
$$\ch(V)=\ch(U)+\ch(W).$$
(ii) If $V$ and $W$ are finite-dimensional $\y$-modules,
$$\ch(V\ot W)=\ch(V)\ch(W). \qed$$
\endproclaim

Part (i) plus Jordan--H\"older means that it suffices to compute the characters
of the irreducible modules $V(P)$ ($P\in\calp$). To state the character
formula, define, for $r\ge 1$, $a\in\Bbb C$,
$$\align
x_a & = e\left(\frac{u-a+1}{u-a}\right)=
e\left(1+\sum_{k=0}^\infty a^ku^{-k-1}\right)\in\Bbb C[\call],\\
y_{r,a}&=1+\sum_{s=1}^r\sum_{t=1}^sx_{a+r-t}^{-1}x_{a+r-t+1}^{-1}\ \ \text{(and
$y_{r,a}=1$ if $r\le 0$)},\\
m_{r,a}(P)&={\roman{max}}\{n\in\Bbb N\mid
\text{$(u-a)^n,(u-a-1)^n,\ldots,(u-a-r+1)^n$ all divide $P(u)$}\}.
\endalign$$
The main result in this section is

\proclaim{Theorem 4.3} For any $P\in\calp$,
$$\ch(V(P))=e(P)\prod_{r\ge 1, a\in\Bbb C}
\left(\frac{y_{r,a}y_{r-2,a+1}}{y_{r-1,a}y_{r-1,a+1}}\right)^{m_{r,a}(P)}.\tag17$$
\endproclaim

\vskip6pt\noindent{\it Remarks} 1. All but finitely many terms in the product
(17) are equal to one, since $m_{r,a}(P)=0$ unless $r\le{\roman{deg}}(P)$ and
$a$ is a root of $P$.

2. It is clear that, if $Q\in\calp$, we have
$$m_{r,a}(PQ)\ge m_{r,a}(P)+m_{r,a}(Q),$$
but strict inequality may occur (e.g. if $P(u)=u-a$ and $Q(u)=u-a-1$, then
$m_{2,a}(PQ)=1$ but $m_{2,a}(P)=m_{2,a}(Q)=0$).

3. The definition of $m_{r,a}$ may be reformulated in terms of the \lq Yangian
derivative\rq:
$$D_{Y}(P)=P(u+1)-P(u).$$
It is clear that
$$m_{r,a}(P)={\roman{max}}\{n\in\Bbb N\mid \text{$(u-a)^n$ divides
$P, D_Y P,\ldots,D_Y^{r-1}P$}\},$$
i.e. that $m_{r,a}(P)$ is the multiplicity of $a$ as a common root of $P, D_Y
P,\ldots,D_Y^{r-1}P$.
\vskip12pt
Before proving 4.2, we note some consequences. Let ${\roman{res}}:\call\to\Bbb
C$ be the homomorphism given by
$${\roman{res}}(f_\bd)=d_0,$$
and let $\call_{\Bbb Z}={\roman{res}}^{-1}(\Bbb Z)$. We also denote by
${\roman{res}}$ the corresponding algebra homomorphisms $\Bbb C[\call]\to\Bbb
C[\Bbb C]$ and $\Bbb C[\call_{\Bbb Z}]\to\Bbb C[\Bbb Z]$. For any
finite-dimensional $\y$-module $V$, we have $\ch(V)\in\Bbb C[\call_{\Bbb Z}]$.
Indeed, this follows from 4.3 if $V$ is irreducible since, for any $P\in\calp$,
$r\ge 1$, $a\in\Bbb C$,
$${\roman{res}}(e(P))=e({\roman{deg}}(P)),\ \
{\roman{res}}(y_{r,a})=\sum_{s=0}^r e(-2s)=z_r\ \text{(say)},$$
and the general case follows from 4.2 (i). Now,
${\roman{res}}(\ch(V))=\ch_{\ung}(V)$, the character of $V$ regarded as an
$\ung$-module. Noting that
$$\sum_{a\in\Bbb C}m_{1,a}(P)={\roman{deg}}(P)$$
and, for $r>1$,
$$\sum_{a\in\Bbb C}m_{r,a}(P)=\text{total number of strings of length $r$ in
$P$}=m_r(P),$$
say, we obtain

\proclaim{Corollary 4.4} For any $P\in\calp$,
$$\ch_{\ung}(V(P))=e({\roman{deg}}(P))z_1^{{\roman{deg}}(P)}\prod_{r=2}^\infty
\left(\frac{z_rz_{r-2}}{z_{r-1}^2}\right)^{m_r(P)}. \qed$$
\endproclaim

Similarly, applying the augmentation homomorphism $\Bbb C[\call]\to\Bbb C$ to
4.3 gives

\proclaim{Corollary 4.5} For any $P\in\calp$,
$${\roman{dim}}(V(P))=2^{{\roman{deg}}(P)}\prod_{r=2}^\infty
\left(\frac{r^2-1}{r^2}\right)^{m_r(P)}. \qed$$
\endproclaim

It follows from 3.9 that, for any $P\in\calp$, $\ch(V(P))$ is an alternating
sum of tensor products of the characters
$$\chi_a=\ch(W_1(a)).$$
The next result makes this explicit. By 3.12,
$$\ch(V(P))=\prod_{r\ge 1,a\in\Bbb C}\ch(W_r(a))^{N_{r,a}(P)},\tag18$$
where $N_{r,a}(P)$ is the number of strings of length $r$ beginning at $a$ in
the canonical decomposition of $S(P)$, so it suffices to consider $P=P_{r,a}$.

\proclaim{Proposition 4.6} For any $r\ge 1$, $a\in\Bbb C$,
$$\ch(W_r(a))=\sum_{s=0}^{[r/2]}(-1)^sA_{r,a}^{(s)},$$
where
$$A_{r,a}^{(s)}=\sum\chi_{a+t_1}\chi_{a+t_2}\ldots\chi_{a+t_{r-2s}},$$
the sum being over those integers $t_1,t_2,\ldots,t_{r-2s}$ such that
$$r>t_1>t_2>\cdots>t_{r-2s}\ge 0\ \ \text{and}\ \ t_j\equiv r-j\ \text{(mod
$2$) for all $j$}.$$
\endproclaim

We shall prove this result after proving the next proposition, which is also
the first step in the proof of 4.3.

\proclaim{Proposition 4.7} For any $r\ge 1$, $a\in\Bbb C$,
$$\ch(W_r(a))=e(P_{r,a})y_{r,a}.$$
\endproclaim
\demo{Proof} From 3.6, we read off that the joint eigenvalue of $w_s\in W_r(a)$
is $\bd_s=\{d_{k,s}\}_{k\in\Bbb N}$, where
$$d_{k,s}=(a+s-1)^ks(r-s+1)-(a+s)^k(s+1)(r-s).$$
This gives
$$f_{\bd_s}=\frac{(u-a+1)(u-a-r)}{(u-a-s+1)(u-a-s)},$$
so
$$\ch(W_r(a))=\sum_{s=0}^r
e\left(\frac{(u-a+1)(u-a-r)}{(u-a-s+1)(u-a-s)}\right).$$
Now,
$$e(P_{r,a})=e\left(\frac{u-a+1}{u-a-r+1}\right),$$
so
$$\ch(W_r(a))=e(P_{r,a})\sum_{s=0}^k
e\left(\frac{u-a-r}{u-a-r+s}\right)
e\left(\frac{u-a-r+1}{u-a-r+s+1}\right)$$
(changing the summation index from $s$ to $r-s$). Since
$$e\left(\frac{u-a-r}{u-a-r+s}\right)=\prod_{t=1}^s
e\left(\frac{u-a-r+t-1}{u-a-r+t}\right)=\prod_{t=1}^s x_{a+r-t+1}^{-1},$$
this gives the stated fromula. \qed\enddemo

\demo{Proof of 4.6} Using 4.7, it is easy to show that
$$\ch(W_{r+2}(a))=\ch(W_{r+1}(a+1))\chi_a-\ch(W_r(a+2)).$$
The formula in 4.6 follows easily from this relation, by using induction on
$r$. \qed\enddemo

If $P\in\calp$, let
$$S(P)=\bigcup_{i=1}^p S_{r_i}(a_i)\tag19$$
be the canonical decomposition of its multiset of roots $S(P)$. By 3.13 and
4.7, we get
$$\ch(V(P))=e(P)\prod_{i=1}^py_{r_i,a_i}.$$
Hence,
$$\ch(V(P))=e(P)\prod_{r\ge 1,a\in\Bbb C}y_{r,a}^{N_{r,a}(P)}.\tag20$$

Now let $n_{r,a}(P)$ be the total number of strings of length $r$ beginning at
$a$ that are contained in $S(P)$. More precisely, relative to the canonical
decomposition (19),
$$n_{r,a}(P)=\sum_{i=1}^p n_i,$$
where
$$n_i=\cases 1&\ \text{if $S_r(a)\subseteq S_{r_i}(a_i)$},\\
0&\ \text{otherwise}.\endcases$$
It is clear that
$$\align
n_{r,a}(P)= &\phantom{+,}N_{r,a}(P)\\
& +N_{r+1,a}(P)+N_{r+1,a-1}(P)\\
& +N_{r+2,a}(P)+N_{r+2,a-1}(P)+N_{r+2,a-2}\\
& +\cdots\\
= & \sum_{s=0}^\infty\sum_{t=0}^sN_{r+s,a-t}(P)\endalign$$
(of course, all but finitely many terms in the double sum are zero). These
equations are easily inverted to express the $N$'s in terms of the $n$'s:
$$N_{r,a}(P)=n_{r,a}(P)-n_{r+1,a}(P)-n_{r+1,a-1}(P)+n_{r+2,a-1}(P).$$
Inserting this into (20), we get
$$\ch(V(P))=e(P)\prod_{r\ge 1,a\in\Bbb C}
\left(\frac{y_{r,a}y_{r-2,a+1}}{y_{r-1,a}y_{r-1,a+1}}\right)^{n_{r,a}(P)}.$$
Thus, 4.3 is a consequence of
\proclaim{Proposition 4.8} For any $m\ge 1$, $a\in\Bbb C$, $P\in\calp$,
$$m_{r,a}(P)=n_{r,a}(P).$$
\endproclaim

To prove 4.8, we need two lemmas.

\proclaim{Lemma 4.9} Let $P\in\calp$. Then, $P$ has a factorisation
$$P=P_1P_2\ldots P_k,$$
for some $k\ge 1$, such that, if $a$ is a root of $P_i$ and $b$ a root of
$P_j$, then $a-b\in\Bbb Z$ if $i=j$ and $a-b\notin\Bbb Z$ if $i\ne j$. \qed
\endproclaim

This is clear.

\proclaim{Lemma 4.10} Let $P\in\calp$ be such that every pair of roots of $P$
differ by an integer. Then, $P$ has a factorisation
$$P=P_1P_2\ldots P_k$$
such that
\vskip6pt\noindent(i) for all $i=1,\ldots,k$, every string in the canonical
decomposition of $S(P_i)$ is a string in the canonical decomposition of $S(P)$;

\noindent (ii) $S(P_i)\cap S(P_j)=\emptyset$ if $i\ne j$;

\noindent (iii) for all $i=1,\ldots,k$, $S(P_i)$ regarded as a set {\rm
without} multiplicities, is a string.
\endproclaim
\demo{Proof} By induction on ${\roman{deg}}(P)$. If  ${\roman{deg}}(P)=0$ or
$1$, there is nothing to prove. Let $S$ be a string of maximal length in the
canonical decomposition of $S(P)$, and let $S_1$ be the union, in the sense of
sets with multiplicities, of all the strings in the canonical decomposition of
$S(P)$ that are contained in $S$. Let $P_1\in\calp$ be the factor of $P$ such
that $S(P_1)=S_1$. By the induction hypothesis,
$$P/P_1=P_2\ldots P_k,$$
where $P_2,\ldots,P_k$ satisfy the conditions of the lemma.

To prove that the factorisation
$$P=P_1P_2\ldots P_k$$
satisfies the conditions of the lemma, we have to prove that, if $j\ne 1$,
\vskip6pt\noindent (a) $S(P_1)\cap S(P_j)=\emptyset$;

\noindent(b) every string in the canonical decomposition of $S(P_1)$ is in
general position with respect to every string in the canonical decomposition
of $S(P_j)$.
\vskip6pt For (a), suppose for a contradiction that there exists $c\in
S(P_1)\cap S(P_j)$. Then, $c\in S\cap T$ for some string $S$ in the canonical
decomposition of $S(P_1)$ and some string $T$ in the canonical decomposition of
$S(P_j)$. But $S$ and $T$ are in general position by construction, so since $S$
has maximal length and $S\cap T\ne\emptyset$, we must have $T\subseteq S$. This
contradicts the definition of $P_1$.

For (b), suppose for a contradiction that $S'$ is a string in the canonical
decomposition of $S(P_1)$, $T'$ a string in the canonical decomposition of
$S(P_j)$, where $j\ne 1$, and that $S'$ and $T'$ are in special position. By
the argument used in the previous paragraph, $S\cap T'=\emptyset$. This gives
two possibilities: either $S$ and $S'$ both begin at some $c\in\Bbb C$ and $T'$
ends at $c-1$, or $S$ and $S'$ both end at some $c'\in\Bbb C$ and $T'$ begins
at $c'+1$. In both cases, $S$ and $T'$ are in special position, a
contradiction.
\qed\enddemo

\demo{Proof of 4.8} By induction on the number of strings in the canonical
decomposition of $S(P)$. The induction begins with
\vskip6pt\noindent{\it Case I: $S(P)$ is a string (multiplicities counted)}.
Then,
$$P(u)=(u-b)(u-b-1)\cdots(u-b-k+1)$$
for some $b\in\Bbb C$, $k\ge 1$. In this case, it is easy to see that
$m_{r,a}(P)$ and $n_{r,a}(P)$ are both equal to $1$ if $r\le k$ and
$a=b,b+1,\ldots,$ or $b+k$, and equal to $0$ otherwise.
\vskip6pt\noindent{\it Case II: $S(P)$ is a string (disregarding
multiplicities)}. Let $S'$ be a string of maximal length in the canonical
decomposition of $S(P)$ (so that $S(P)=S'$ disregarding multiplicities), and
let $S''$ be the union of the other strings in the canonical decomposition of
$S(P)$. Let $P=P'P''$ be the corresponding factorisation of $P$. It is clear
from the definition of $n_{r,a}$ that
$$n_{r,a}(P)=n_{r,a}(P')+n_{r,a}(P''),$$
and since $P'$ has no repeated roots,
$$m_{r,a}(P)=m_{r,a}(P')+m_{r,a}(P''),$$
so the result follows by induction.
\vskip6pt\noindent{\it Case III: Any two roots of $P$ differ by an integer}. We
have a factorisation
$$P=P_1P_2\ldots P_k,$$
where the $P_i$ satsify the condition in 4.10. Since $S(P_i)\cap
S(P_j)=\emptyset$ if $i\ne j$, it is clear that, for each $r\ge 1$, $a\in\Bbb
C$, $n_{r,a}(P_i)>0$ for at most one $i$, say $i=1$ without loss of generality,
and that $n_{r,a}(P)=n_{r,a}(P_1)$. Hence,
$$n_{r,a}(P)=\sum_{i=1}^k n_{r,a}(P_i).$$

On the other hand, conditions (i) and (ii) in 4.10 imply that, if $a\in
S(P_i)$, $b\in S(P_j)$, and $i\ne j$, then $|a-b|\ge 2$. We claim that this
implies that
$$m_{r,a}(P)=\sum_{i=1}^k m_{r,a}(P_i),\tag21$$
so that Case III follows from Case II.

Suppose that $(u-a)^m, (u-a-1)^m,\ldots,(u-a-r+1)^m$ all divide $P$. Since
$S(P_i)\cap S(P_j)=\emptyset$ if $i\ne j$, $(u-a)^m$ divides $P_j$ for some
$j$. Similarly, $(u-a-1)^m$ divides $P_i$ for some $i$, and we must have $i=j$
otherwise a root of $P_j$ would differ from a root of $P_i$ by less than $2$.
Continuing in this way, we see that $(u-a)^m$, $(u-a-1)^m,\ldots,(u-a-r+1)^m$
all divide $P_j$ and divide no other $P_i$. This proves that $m_{r,a}(P_i)=0$
if $i\ne j$, and $m_{r,a}(P)\le m_{r,a}(P_j)$. Hence,
$$m_{r,a}(P)\le\sum_{i=1}^km_{r,a}(P_i).$$
The opposite inequality is obvious (see Remark 2 following 4.3), so (21) is
proved.
\vskip6pt\noindent{\it Case IV: General case}. By 4.9, we have a factorisation
$$P=P_1P_2\ldots P_k,$$
where each $P_i$ satisfies the hypotheses of Case III and, if $i\ne j$, each
root of $P_i$ differs by a non-integer from each root of $P_j$. It is now clear
that
$$m_{r,a}(P)=\sum_{i=1}^k m_{r,a}(P_i),\ \ \
n_{r,a}(P)=\sum_{i=1}^k n_{r,a}(P_i),$$
(at most one term in each sum being non-zero), so the result follows from Case
III. \qed\enddemo
\vskip36pt\centerline{\bf References}
\vskip12pt\noindent 1.  D. Bernard and A. LeClair, Quantum group symmetries and
non-local currents in 2D QFT, Commun. Math. Phys. {\bf 142} (1991), 99--138.

\noindent 2. V. Chari and A. N. Pressley, New unitary representations of loop
groups, Math. Ann. {\bf 275} (1986), 87--104.

\noindent 3.  V. Chari and A. N. Pressley, Yangians and R-matrices, L'Enseign.
Math. {\bf 36} (1990), 267--302.

\noindent 4. V. Chari and A. N. Pressley, Fundamental representations of
Yangians and singularities of R-matrices, J. reine angew. Math. {\bf 417}
(1991), 87--128.

\noindent 5. V. Chari and A. N. Pressley, {\it A Guide to Quantum Groups},
Cambridge University Press, Cambridge, 1994.

\noindent 6. V. Chari and A. N. Pressley, Minimal affinizations of
representations of quantum groups: the irregular case, preprint, 1995.

\noindent 7. V. Drinfel'd, Hopf algebras and the quantum Yang--Baxter equation,
Soviet Math. Dokl. {\bf 32} (1985), 254--258.

\noindent 8. V. Drinfel'd, A new realization of Yangians and quantum affine
algebras, Soviet Math. Dokl. {\bf 36} (1988), 212--216.

\noindent 9. H. Knight, Spectrum of finite-dimensional representations of
Yangians, preprint, Yale University, 1994.
\vskip36pt
{\eightpoint{
$$\matrix\format\l&\l&\l&\l\\
\phantom{.} & {\text{Vyjayanthi Chari}}\phantom{xxxxxxxxxxxxx} & {\text{Andrew\
Pressley}}\\
\phantom{.}&{\text{Department of Mathematics}}\phantom{xxxxxxxxxxxxx} &
{\text{Department of Mathematics}}\\
\phantom{.}&{\text{University of California}}\phantom{xxxxxxxxxxxxx} &
{\text{King's College}}\\
\phantom{.}&{\text{Riverside}}\phantom{xxxxxxxxxxxxx} & {\text{Strand}}\\
\phantom{.}&{\text{CA 92521}}\phantom{xxxxxxxxxxxxx} & {\text{London WC2R
2LS}}\\
\phantom{.}&{\text{USA}}\phantom{xxxxxxxxxxxxx} & {\roman{UK}}\\
&{\text{email: chari\@ucrmath.ucr.edu}}\phantom{xxxxxxxxxxxxx} &
{\text{email:anp\@uk.ac.kcl.mth}}
\endmatrix$$
}}

\enddocument